# Statistical inference under order restrictions on both rows and columns of a matrix, with an application in toxicology

Eric Teoh[*,1], Abraham Nyska[2], Uri Wormser[3] and Shyamal D. Peddada[†,4]

*Insurance Institute for Highway Safety, Tel Aviv University, The Hebrew University of Jerusalem and NIEHS*

**Abstract:** We present a general methodology for performing statistical inference on the components of a real-valued matrix parameter for which rows and columns are subject to order restrictions. The proposed estimation procedure is based on an iterative algorithm developed by Dykstra and Robertson (1982) for simple order restriction on rows and columns of a matrix. For any order restrictions on rows and columns of a matrix, sufficient conditions are derived for the algorithm to converge in a single application of row and column operations. The new algorithm is applicable to a broad collection of order restrictions. In practice, it is easy to design a study such that the sufficient conditions derived in this paper are satisfied. For instance, the sufficient conditions are satisfied in a balanced design. Using the estimation procedure developed in this article, a bootstrap test for order restrictions on rows and columns of a matrix is proposed. Computer simulations for ordinal data were performed to compare the proposed test with some existing test procedures in terms of size and power. The new methodology is illustrated by applying it to a set of ordinal data obtained from a toxicological study.

## 1. Introduction

In many applications the parameter of interest $\theta$ can be expressed as elements of a real-valued $I \times J$ matrix such that the elements within rows and/or columns of the matrix are subject to inequality restrictions called order restrictions. Researchers are often interested in drawing statistical inferences on such parameters subject to a variety of order restrictions. For a parameter vector $\eta = (\eta_1, \eta_2, \ldots, \eta_p)'$, some commonly encountered order restrictions are; *simple order*, where $\eta_1 \leq \eta_2 \leq \cdots \leq$

[*]This work was done while the author was a summer intern at the National Institute of Environmental Health Sciences (NIEHS) and a graduate student in the Department of Biostatistics, University of North Carolina, Chapel Hill, NC, USA.

[†]Supported by the Intramural Research Program of the NIH, National Institute of Environmental Health Sciences.

[1]Insurance Institute for Highway Safety, Arlington, Virginia, USA, e-mail: eteoh@iihs.org
[2]Toxicologic Pathologist, Tel Aviv University, Tel Aviv, Israel, e-mail: anyska@bezeqint.net
[3]Department of Pharmacology, School of Pharmacy, Faculty of Medicine, Institute of Life Sciences, The Hebrew University of Jerusalem, Jerusalem, Israel, e-mail: wormser@cc.huji.ac.il
[4]Biostatistics Branch, NIEHS, Research Triangle Park, NC, USA, e-mail: peddada@niehs.nih.gov

Table 1
*For each response variable, the data structure of the experiment in Wormser et al. [26]*

| Genotype | Sample size | Level of skin injury | | | | |
|---|---|---|---|---|---|---|
| | | Unremarkable | Minimal | Mild | Moderate | Marked |
| Cox-2 deficient | 10 | $n_{1,1}$ | $n_{1,2}$ | $n_{1,3}$ | $n_{1,4}$ | $n_{1,5}$ |
| Wild type | 10 | $n_{2,1}$ | $n_{2,2}$ | $n_{2,3}$ | $n_{2,4}$ | $n_{2,5}$ |
| Cox-1 deficient | 10 | $n_{3,1}$ | $n_{3,2}$ | $n_{3,3}$ | $n_{3,4}$ | $n_{3,5}$ |

$\eta_p$, *umbrella order*, where $\eta_1 \leq \eta_2 \leq \cdots \eta_i \geq \eta_{i+1} \geq \cdots \geq \eta_p$ and *simple tree order*, where $\eta_1 \leq \eta_i$, for all $2 \leq i \leq p$.

Recently Wormser et al. [26] conducted an experiment to evaluate the differences among three different genotypes of mice, namely, the wild type (WT), the cyclooxygenase-1 deficient (COX-1-d) and the cyclooxygenase-2 deficient (COX-2-d) mice, when they were exposed to sulfur mustard (also known as mustard gas). Depending upon the severity of injury to skin, each animal was categorized into one of five ordered categories, namely, "unremarkable", "minimal", "mild", "moderate", and "marked" (details in Section 5). In this experiment, a sample of $n = 10$ animals from each genotype was exposed to sulfur mustard. Let $n_{i,j}$, $i = 1, 2, 3$, $j = 1, 2, \ldots, 5$, denote the number of animals in the $i^{th}$ genotype that belong to $j^{th}$ response category, with $E(n_{i,j}) = n\pi_{i,j}$. Then the parameters of interest are the cumulative probabilities $\theta_{i,j} = \sum_{k=1}^{j} \pi_{i,k}$, $i = 1, 2, 3$, $j = 1, 2, 3, 4$. Note that $\theta_{i,5} = 1$, $i = 1, 2, 3$. Table 1 summarizes the type of data obtained in Wormser et al. [26]. Clearly, the rows of $\theta$ satisfy a simple order as they represent cumulative probabilities. According to Wormser et al. [26], COX-2 deficiency has a protective effect against inflammatory processes while COX-1 deficiency has a negative effect. Consequently, each column of $\theta$ is also subject to simple order restriction. Thus in this example the rows as well as columns of $\theta$ are subject to simple order restriction.

The above type of matrix order restrictions commonly arise in a variety contexts such as the analysis of ordinal data (cf. Agresti and Coull [1], Grove [9], Nair [16] and Wang [25]), survival analysis (Praestgaard [19]), Phase I clinical trials involving "cocktail" of treatments (Conaway et al. [6]) and analysis of time-course and dose-response gene expression microarray studies, etc.

For a given parametric family with matrix valued parameter $\theta \in \Theta \subset \mathcal{R}^{I \times J}$, where $\Theta$ is the parameter space defined by order restrictions on the rows and/or columns of $\theta$, one may estimate $\theta$ using restricted maximum likelihood estimators (RMLEs) and test alternative hypotheses using the likelihood ratio tests and their modifications. Such methods have been well studied in the literature (cf. Barlow et al. [2], Robertson et al. [21] and Silvapulle and Sen [23]). However, as seen from Hwang and Peddada [11] and Lee [12], the RMLE is not always efficient relative to the unrestricted maximum likelihood estimator (UMLE). Also, the likelihood ratio tests in the present context may not be computationally simple and the asymptotic distribution under the null hypothesis may involve nuisance parameters (cf. Franck [8], Grove [9], Robertson and Wright [20] and Wang [25]). Silvapulle [22] provided an interesting explanation for why sometimes the RMLEs and the likelihood ratio tests provide counter-intuitive results.

In view of the general concerns regarding RMLE and the likelihood ratio tests, in Section 2 we introduce a computationally straightforward methodology for estimating a matrix valued parameter $\theta$ when the rows and/or columns are subject to order restrictions. The proposed estimation procedure is based on the iterative procedure of Dykstra and Robertson [7] and uses the method of Hwang and Peddada [11] for general order restrictions. If the elements of the random matrix are independently and normally distributed and if rows as well as columns of the matrix are subject to simple order restriction then the proposed procedure is identical to the estimator given in Dykstra and Robertson [7], but the two procedures may differ for other order restrictions. Using the proposed point estimators, in Section 3 we introduce a Kolmogorov-Smirnov type test statistic and a bootstrap-based methodology for determining significance. The performance of the proposed test procedure is evaluated using computer simulations in Section 4 and in Section 5 it is illustrated by analyzing the data in Wormser et al. [26] mentioned above. Concluding remarks are provided in Section 6.

## 2. Estimation of parameters subject to order restrictions

### 2.1. Notations and a brief review

Throughout this paper $\mathcal{R}^p$ denotes the vector space of $p \times 1$ real vectors and $\mathcal{R}^{I \times J}$ denotes the vector space of $I \times J$ real matrices. Two components $\theta_{i,j}$ and $\theta_{r,s}$ of $\theta \in \Theta \subset \mathcal{R}^{I \times J}$ are said to be *linked* if the inequality between them is known *a priori*. A parameter is said to be a *nodal* parameter if it is linked to all $IJ$ components of $\theta \in \Theta \subset \mathcal{R}^{I \times J}$. A subset of parameters $\mathcal{M}$, formed by the components of $\theta \in \Theta \subset \mathcal{R}^{I \times J}$, is a *linked subgraph* if all parameters in $\mathcal{M}$ are linked, with at least one strict inequality. Note that every linked subgraph represents a simple order-restriction and conversely, every simple order is a linked subgraph. A linked subgraph $\mathcal{M}$ is said to be *maximally linked* if for any linked subgraph $\mathcal{N}$, $\mathcal{M} \subset \mathcal{N} \implies \mathcal{M} = \mathcal{N}$.

If a subset of parameters simultaneously satisfies two linked subgraphs $\mathcal{M}$ and $\mathcal{N}$, then we use the notation $\mathcal{M}\Lambda\mathcal{N}$ to describe the subset. As noted in Peddada et al. [18], any order-restriction between parameters can be expressed in terms of a collection of maximally linked subgraphs. Similarly, every parameter $\theta_{i,j}$ appears in a finite collection of maximally linked subgraphs as a nodal parameter (within each subgraph). Such maximally linked subgraphs are said to be *associated* with $\theta_{i,j}$.

If the rows of $\theta$ are subject to order restriction $\mathcal{R} \subset \mathcal{R}^I$ and the columns are subject to order restriction $\mathcal{C} \subset \mathcal{R}^J$ then we shall denote the joint order restriction by $\mathcal{R}\Lambda\mathcal{C} \subset \mathcal{R}^{I \times J}$.

The notation $\mathcal{R}\Lambda\mathcal{R}^J$ would indicate order restrictions $\mathcal{R}$ on the rows and no order restrictions on the columns of $\theta$. Similarly, the notation $\mathcal{R}^I\Lambda\mathcal{C}$ would indicate order restrictions $\mathcal{C}$ on the columns and no order restrictions on the rows of $\theta$.

It is important to emphasize that in this article we only consider linked subgraphs which are subsets of rows or which are subsets of columns of the parameter $\theta$. Thus we are not considering inequalities between arbitrary linear or non-linear functions of the elements of $\theta$.

**Example 2.1** (Umbrella order restriction on the rows of $\theta$). Suppose the elements of each row of $\theta$ are subject to an umbrella order with peak in the $s^{th}$ column. That is, the components of the $i^{th}$ row of $\theta$ satisfy the inequalities $\theta_{i,1} \leq \theta_{i,2} \leq \cdots \leq \theta_{i,s} \geq \theta_{i,s+1} \geq \cdots \geq \theta_{i,J}$, $i = 1, 2, \ldots, I$. Then in this case the parameters of the $i^{th}$ row can be expressed in terms of two maximally linked subgraphs, namely, $\mathcal{M}_{i,1:s} = \{(\theta_{i,1}, \theta_{i,2}, \ldots \theta_{i,s}) \mid \theta_{i,1} \leq \theta_{i,2} \leq \cdots \leq \theta_{i,s}\}$ and $\mathcal{M}_{i,s:J} = \{(\theta_{i,s}, \theta_{i,s+1}, \ldots, \theta_{i,J}) \mid \theta_{i,s} \geq \theta_{i,s+1} \geq \cdots \geq \theta_{i,J}\}$. Thus, the subset of parameters in the $i^{th}$ row of $\theta$, $i = 1, 2, \ldots, I$, can be expressed as $\mathcal{M}_{i,1:s}\Lambda\mathcal{M}_{i,s:J}$. Further, for $r \leq s$, the maximally linked subgraph associated with $\theta_{i,r}$ in the $i^{th}$ row, $i = 1, 2, \ldots, I$, is $\mathcal{M}_{i,1:s}$.

For a vector $x = (x_1, x_2, \ldots, x_p)' \in \mathcal{R}^p$, a *simple order operator* $\mathfrak{C}_w^\mathcal{S} : \mathcal{R}^p \to \mathcal{S}$ is an orthogonal projection operator onto the simple order cone $\mathcal{S} = \{x \in \mathcal{R}^p : x_1 \leq x_2 \leq \cdots \leq x_p\}$ where $w$ is a vector of positive weights and the $i^{th}$ component of $\mathfrak{C}_w^\mathcal{S}(x)$ is given by

$$(2.1) \qquad (\mathfrak{C}_w^\mathcal{S}(x))_i = \min_{s \geq i} \max_{t \leq i} \frac{\sum_{j=t}^s w_j x_j}{\sum_{j=t}^s w_j}.$$

We shall use the terminology *simple order function* to describe the min-max formula used in (2.1).

**Remark 2.1.** For a fixed weight vector $w$, $\mathfrak{C}_w^\mathcal{S}$ is a monotonic operator in the sense that if $x \leq y$ then $\mathfrak{C}_w^\mathcal{S}(x) \leq \mathfrak{C}_w^\mathcal{S}(y)$, where the inequality is componentwise.

When estimating a parameter $\eta = (\eta_1, \eta_2, \ldots, \eta_p)'$ subject to arbitrary order restrictions $\mathcal{C} \subset \mathcal{R}^p$ with at least one nodal parameter, Hwang and Peddada [11] used the simple order operator $\mathfrak{C}_w^\mathcal{S}$, with a suitable weight vector $w$, for estimating the nodal parameters. Typically the weights are proportional to the precision (or sometimes the sample size) of UMLE. Once a parameter is estimated, then the corresponding UMLE is replaced by the new restricted estimator and is assigned arbitrarily large weight $B$, $B \to \infty$, in all subsequent calculations. To estimate a non-nodal parameter, identify the collection of all maximally linked subgraphs associated with that non-nodal parameter and apply the simple order operator $\mathfrak{C}_w^\mathcal{S}$ on the vector corresponding to the subgraph with a suitable weight vector $w$. Suitable modifications were proposed for graphs with no nodal parameters.

**Example 2.2** (Umbrella order restriction). Suppose $\eta$ is a parameter satisfying the order restriction $\eta_1 \leq \eta_2 \leq \eta_3 \geq \eta_4 \geq \eta_5$ with UMLE $\hat{\eta}$. In this case the only nodal parameter is $\eta_3$. Suppose $w = (w_1, w_2, w_3, w_4, w_5)'$ is the weight vector associated with $\hat{\eta}$. According to Hwang and Peddada [11] the estimation procedure begins with the nodal parameter $\eta_3$. For $i = 1, 2$, denote $w_{(i)} = w_i, \hat{\eta}_{(i)} = \hat{\eta}_i$ and let $w_{(3)} = w_5, w_{(4)} = w_4, w_{(5)} = w_3$ and $\hat{\eta}_{(3)} = \hat{\eta}_5, \hat{\eta}_{(4)} = \hat{\eta}_4, \hat{\eta}_{(5)} = \hat{\eta}_3$. Then $\eta_3$ may be estimated by the following simple order formula

$$\hat{\hat{\eta}}_3 = \max_{i \leq 5} \frac{\sum_{j=i}^5 w_{(j)} \hat{\eta}_{(j)}}{\sum_{j=i}^5 w_{(j)}}.$$

Next, the non-nodal parameters $\eta_1, \eta_2, \eta_4$ and $\eta_5$ are estimated using the maximally linked subgraphs $\eta_1 \leq \eta_2 \leq \eta_3$ and $\eta_3 \geq \eta_4 \geq \eta_5$, respectively. The estimators for the non-nodal parameters $\eta_1, \eta_2$ are simplified as follows:

$$\hat{\hat{\eta}}_1 = \min\{\hat{\eta}_1, \frac{w_1 \hat{\eta}_1 + w_2 \hat{\eta}_2}{w_1 + w_2}, \hat{\hat{\eta}}_3\}, \qquad \hat{\hat{\eta}}_2 = \min\{\max(\frac{w_1 \hat{\eta}_1 + w_2 \hat{\eta}_2}{w_1 + w_2}, \hat{\eta}_2), \hat{\hat{\eta}}_3\}.$$

In a similar manner $\eta_4$ and $\eta_5$ are estimated.

**Remark 2.2.** For a given $w$, $\mathfrak{C}_w^\mathcal{C}$ is a function of several simple order functions and therefore $\mathfrak{C}_w^\mathcal{C}$ is a monotonic operator. That is, for all $x \leq y$, $\mathfrak{C}_w^\mathcal{C}(x) \leq \mathfrak{C}_w^\mathcal{C}(y)$, where the inequalities are componentwise.

### 2.2. Estimation of matrix valued parameters

We extend the notations from the previous section to matrix valued parameters as follows. For a matrix $X \in \mathcal{R}^{I \times J}$ and a weight matrix $W$, we denote the *matrix* *simple order column operator* by $\mathfrak{C}_W^{\mathcal{S}}$. Each column $x$ of $X$ is orthogonally projected onto the simple order cone $\mathcal{S}$ using $\mathfrak{C}_w^{\mathcal{S}}$, where $w$ is a suitable column vector of $W$. Similarly, the *matrix simple order row operator* is denoted by $\mathfrak{R}_W^{\mathcal{S}}$, which, with rows of $W$, projects each row vector of matrix $X \in \mathcal{R}^{I \times J}$ orthogonally onto the simple order cone $\mathcal{S}$. Analogously, for arbitrary order restrictions $\mathcal{D}$, we define matrix column and row operators using Hwang and Peddada methodology by $\mathfrak{C}_W^{\mathcal{D}}$ and $\mathfrak{R}_W^{\mathcal{D}}$, respectively.

We now describe the proposed algorithm for estimating $\theta \in \Theta = \mathcal{R}\Lambda\mathcal{C} \subset \mathcal{R}^{I \times J}$ using an unrestricted point estimator $\hat{\theta}$. We use the weight matrix $W_R$ when operating on the rows of $\hat{\theta}$ and weight matrix $W_C$ when operating on the columns of $\hat{\theta}$.

*Step 1 (An unrestricted estimator)*:
Obtain an unrestricted estimator $\hat{\theta}$ for an $I \times J$ matrix parameter $\theta$. Note that in most situations a user may prefer to start with the unrestricted maximum likelihood estimator (UMLE), although it is not required.

*Step 2 (Estimation under order restrictions on the columns of $\theta$)*:
Apply the procedure of Hwang and Peddada [11] on each column of $\hat{\theta}$ to obtain estimates for $\theta$ under the order-restriction on the columns of $\theta$. That is, apply the operator $\mathfrak{C}_{W_C}^{\mathcal{C}}$ on the columns of $\hat{\theta}$ and denote the resulting estimator by $\mathfrak{C}_{W_C}^{\mathcal{C}}(\hat{\theta})$. Note that the elements of $\mathfrak{C}_{W_C}^{\mathcal{C}}(\hat{\theta})$ may not satisfy the order-restriction on the rows of $\theta$.

*Step 3 (Estimation under order restrictions on the rows of $\theta$)*:
Apply the operator $\mathfrak{R}_{W_R}^{\mathcal{R}}$ on the rows of $\mathfrak{C}_{W_C}^{\mathcal{C}}(\hat{\theta})$ and denote the resulting estimator by $(\mathfrak{R}_{W_R}^{\mathcal{R}} \circ \mathfrak{C}_{W_C}^{\mathcal{C}})(\hat{\theta})$.

*Step 4 (Iterate to convergence)*:
Repeat Steps 2 and 3 to obtain the $q^{th}$ iterate, $(\mathfrak{R}_{W_R}^{\mathcal{R}} \circ \mathfrak{C}_{W_C}^{\mathcal{C}})^q(\hat{\theta})$. Stop when some reasonable convergence criterion is reached. Denote $\lim_{q \to \infty} (\mathfrak{R}_{W_R}^{\mathcal{R}} \circ \mathfrak{C}_{W_C}^{\mathcal{C}})^q(\hat{\theta})$ as $\tilde{\theta}_1$.

*Step 5 (Final estimate)*:
Since $\mathfrak{R}_{W_R}^{\mathcal{R}}$ and $\mathfrak{C}_{W_C}^{\mathcal{C}}$ do not necessarily commute for all order restrictions and weights used in the calculated weighted averages, repeat the process beginning with within-row order restrictions followed by within-column order restrictions. That is, compute $\tilde{\theta}_2 = \lim_{q \to \infty} (\mathfrak{C}_{W_C}^{\mathcal{C}} \circ \mathfrak{R}_{W_R}^{\mathcal{R}})^q(\hat{\theta})$. The final estimate of $\theta \in \Theta$ is taken to be

$$(2.2) \qquad \hat{\hat{\theta}} \equiv \frac{1}{2}(\tilde{\theta}_1 + \tilde{\theta}_2).$$

**Remark 2.3.** In general $\tilde{\theta}_1 \neq \tilde{\theta}_2$. To illustrate this, consider the very special case where $\theta$ is a $2 \times 2$ matrix with rows and columns both subject to the simple order restriction $\theta_{1,j} \leq \theta_{2,j}$, $\theta_{i,1} \leq \theta_{i,2}$, $i = 1,2$, and $j = 1,2$. Consider the extreme case where $\hat{\theta}$ is a symmetric matrix with $\hat{\theta}_{1,1} \leq \hat{\theta}_{1,2} = \hat{\theta}_{2,1} \geq \hat{\theta}_{2,2}$. Further, suppose $W_C = W_R = 11'$. Thus we have a perfectly "symmetric" problem where we may expect $\tilde{\theta}_1 = \tilde{\theta}_2$. However, even in this rather seemingly obvious situation $\tilde{\theta}_1 \neq \tilde{\theta}_2$, but $\tilde{\theta}_1 = \tilde{\theta}_2'$. Hence the composition $(\mathfrak{R}_{W_R}^{\mathcal{R}} \circ \mathfrak{C}_{W_C}^{\mathcal{C}})$ is not commutative. For this reason, we need to invoke Step 5 in all situations.

Before we discuss the convergence of the above algorithm in Steps 4 and 5, we consider the following example which may be instructive.

**Example 2.3.** Consider a clinical trial comparing 3 new treatments with an existing treatment using 4 dose groups per treatment group. For the $j^{th}$ dose of the $i^{th}$ treatment, $i, j = 1, 2, \ldots, 4$, let $\hat{\theta}_{i,j}$ denote the sample mean response based on $n_{i,j} = n$ observations. Further, for simplicity of illustration, let $\hat{\theta}_{i,j} \sim^{indep} N(\theta_{i,j}, c)$, $i, j = 1, 2, 3, 4$. Without loss of generality, let the first row of $\theta$ correspond to the existing treatment. Then the order restrictions of interest are $\theta_{1,j} \leq \theta_{i,j}$, $i \geq 2$, $j \geq 1$, i.e. simple tree order within each column of $\theta$, and $\theta_{i,j_1} \leq \theta_{i,j_2}$, $1 \leq j_1 \leq j_2 \leq 4$, $i \geq 1$, i.e. simple order within each row of $\theta$. We choose $W_C = W_R = 11'$, where $1 = (1, 1, 1, 1)'$.

We begin with Step 2 of the algorithm. Let $Y = \mathfrak{C}^{\mathcal{C}}_{W_C}(\hat{\theta})$. Thus columns of $Y$ satisfy simple tree order. That is,

$$Y_{i,j} \geq Y_{1,j}, \forall i = 2, 3, 4, \quad j = 1, 2, 3, 4.$$

Now applying Step 3 on $Y = \mathfrak{C}^{\mathcal{C}}_{W_C}(\hat{\theta})$ we obtain $Z = \mathfrak{R}^{\mathcal{R}}_{W_R}(Y)$. The rows of $Z$ satisfy a simple order. That is,

$$Z_{i,j_1} \geq Z_{i,j_2}, \forall\, j_1, j_2 = 1, 2, 3, 4,\ j_1 > j_2,\ i = 1, 2, 3, 4.$$

We now demonstrate that the columns of $Z$ would also satisfy the simple tree order restriction imposed on $\theta$. That is, for any column $j$ we need to demonstrate that $Z_{i,j} \geq Z_{1,j}$, for all $i = 2, 3, 4$. Note that for each $j = 1, 2, 3, 4$,

$$Z_{i,j} = \min_{t \geq j} \max_{s \leq j} \frac{\sum_{k=s}^{t} Y_{i,k}}{t - s + 1}, \forall i = 1, 2, 3, 4.$$

Since $Y_{i,k} \geq Y_{1,k}, \forall i = 2, 3, 4,\ k = 1, 2, 3, 4$ and since the above simple order function is a monotonic function it follows that for all $i = 2, 3, 4$ and $j = 1, 2, 3, 4$,

$$Z_{i,j} = \min_{t \geq j} \max_{s \leq j} \frac{\sum_{k=s}^{t} Y_{i,k}}{t - s + 1} \geq \min_{t \geq j} \max_{s \leq j} \frac{\sum_{k=s}^{t} Y_{1,k}}{t - s + 1} = Z_{1,j}.$$

Thus $\tilde{\theta}_1 = (\mathfrak{R}^{\mathcal{R}}_{W_R} \circ \mathfrak{C}^{\mathcal{C}}_{W_C})(\hat{\theta})$. Similarly, it can be demonstrated that $\tilde{\theta}_2 = (\mathfrak{C}^{\mathcal{C}}_{W_C} \circ \mathfrak{R}^{\mathcal{R}}_{W_R})(\hat{\theta})$. Thus in this example the algorithm converges after one application of column and row operations and $q = 1$.

As will be demonstrated formally in the following theorem, one of the reasons for the convergence observed in the above example is that $W_C = W_R = 11'$, where $1$ is a column vector of 1's of suitable length. Or more generally, $W_C$ and $W_R$ are each of rank 1. In many applications, researchers use a balanced design for all dose and treatment combinations. In such situations it is appropriate to take $W_C = W_R = 11'$.

We now discuss convergence of Steps 4 and 5 in the above algorithm in the following theorem.

**Theorem 2.1.** *Suppose $\theta \in \Theta = \mathcal{R}\Lambda\mathcal{C} \subset \mathcal{R}^{I \times J}$, with every row subject to the same order restriction and every column subject to the same order restriction. However, the order restriction on a row need not be same as the order restriction on a column. Further, suppose that the weight matrices $W_R$ and $W_C$ are each of rank 1. Then for all $X \in \mathcal{R}^{I \times J}$,*

*(a)* $\mathfrak{C}^{\mathcal{C}}_{W_C} \circ (\mathfrak{R}^{\mathcal{R}}_{W_R} \circ \mathfrak{C}^{\mathcal{C}}_{W_C})(X) = (\mathfrak{R}^{\mathcal{R}}_{W_R} \circ \mathfrak{C}^{\mathcal{C}}_{W_C})(X) \in \Theta.$

*(b)* $\mathfrak{R}^{\mathcal{R}}_{W_R} \circ (\mathfrak{C}^{\mathcal{C}}_{W_C} \circ \mathfrak{R}^{\mathcal{R}}_{W_R})(X) = (\mathfrak{C}^{\mathcal{C}}_{W_C} \circ \mathfrak{R}^{\mathcal{R}}_{W_R})(X) \in \Theta.$

*Thus it takes one column operation and one row operation for the algorithm described in Steps 4 and 5 to converge.*

*Proof.* We prove the theorem for (a) since the proof of (b) follows similarly. The main underlying idea of the proof is that, as stated in Remark 2.2, the row and column operators $\mathfrak{R}^{\mathcal{R}}_{W_R}$ and $\mathfrak{C}^{\mathcal{C}}_{W_C}$ are monotonic and the simple order function used in these operators is a function of suitable weighted averages of the elements of $\hat{\theta}$.

Let $Y = \mathfrak{C}^{\mathcal{C}}_{W_C}(X)$. Thus, the columns of $Y$ satisfy the order restriction $\mathcal{C}$. That is, for each $j = 1, 2, \ldots, J$, and any two linked parameters $\theta_{i_1,j}$ and $\theta_{i_2,j}$, with $\theta_{i_1,j} \leq \theta_{i_2,j}$, we have $Y_{i_1,j} \leq Y_{i_2,j}$. Let $Z = \mathfrak{R}^{\mathcal{R}}_{W_R}(Y) = (\mathfrak{R}^{\mathcal{R}}_{W_R} \circ \mathfrak{C}^{\mathcal{C}}_{W_C})(X)$. Thus, the rows of $Z$ satisfy the order restriction $\mathcal{R}$.

For each $j = 1, 2, \ldots, J$, for any two linked parameters $\theta_{i_1,j}$ and $\theta_{i_2,j}$, with $\theta_{i_1,j} \leq \theta_{i_2,j}$, we demonstrate that $Z_{i_1,j} \leq Z_{i_2,j}$.

Recall from Remark 2.2 that the operators $\mathfrak{C}^{\mathcal{C}}_{W_C}$ and $\mathfrak{R}^{\mathcal{R}}_{W_R}$ are functions of several simple order functions of suitable weighted averages of the components of $\hat{\theta}$. For a given row $i_1$, let $\mathcal{L}$ denote the set of all subsets of $\mathcal{J} = \{1, 2, 3, \ldots, J\}$ such that $\hat{\theta}$ with column indices in these sets are used in the construction of $Z_{i_1,j}$. Further, since the order restrictions in every column is the same, the same set of column indices are used in the construction of $Z_{i_2,j}$. Thus $Z_{i_1,j}$ and $Z_{i_2,j}$ are functions of $\sum_{k \in K}(W_R)_{i_1,k} Y_{i_1,k} / \sum_{k \in K}(W_R)_{i_1,k}$ and $\sum_{k \in K}(W_R)_{i_2,k} Y_{i_2,k} / \sum_{k \in K}(W_R)_{i_2,k}$, respectively, where $K \subset \mathcal{L}$.

Since $\mathfrak{R}^{\mathcal{R}}_{W_R}$ is a monotonic operator, it is sufficient to prove that

$$(2.3) \quad \frac{\sum_{k \in K}(W_R)_{i_1,k} Y_{i_1,k}}{\sum_{k \in K}(W_R)_{i_1,k}} \leq \frac{\sum_{k \in K}(W_R)_{i_2,k} Y_{i_2,k}}{\sum_{k \in K}(W_R)_{i_2,k}}, \text{ for all } K \in \mathcal{L}.$$

Since $W_R$ is of rank 1, we can therefore express $(W_R)_{i_1,k} = \alpha_{i_2}(W_R)_{i_2,k}$. Therefore

$$(2.4) \quad \frac{\sum_{k \in K}(W_R)_{i_1,k} Y_{i_1,k}}{\sum_{k \in K}(W_R)_{i_1,k}} = \frac{\sum_{k \in K} \alpha_{i_2}(W_R)_{i_2,k} Y_{i_1,k}}{\sum_{k \in K} \alpha_{i_2}(W_R)_{i_2,k}} = \frac{\sum_{k \in K}(W_R)_{i_2,k} Y_{i_1,k}}{\sum_{k \in K}(W_R)_{i_2,k}}.$$

But since $Y_{i_1,j} \leq Y_{i_2,j}$ for every $j = 1, 2, \ldots, J$, (2.4) is bounded by

$$\frac{\sum_{k \in K}(W_R)_{i_2,k} Y_{i_2,k}}{\sum_{k \in K}(W_R)_{i_2,k}}.$$

Thus, by the monotonicity of the operator $\mathfrak{R}^{\mathcal{R}}_{W_R}$, $Z_{i_1,j} \leq Z_{i_2,j}$. Hence $(\mathfrak{R}^{\mathcal{R}}_{W_R} \circ \mathfrak{C}^{\mathcal{C}}_{W_C})(X) \in \Theta$ and therefore $\mathfrak{C}^{\mathcal{C}}_{W_C}$ is a left identity of $(\mathfrak{R}^{\mathcal{R}}_{W_R} \circ \mathfrak{C}^{\mathcal{C}}_{W_C})(X)$.

Hence the proof of the theorem. □

**Remark 2.4.** An important consequence of the above theorem is that $W_C$ and $W_R$ need not be limited to the matrix $11'$, but $I \times J$ weight matrices of the type $n_{ij} = n_j$, $i = 1, 2, \ldots, I$, $j = 1, 2, \ldots, J$, can be considered. For example, in a treatment by dose-response study, it is not required to have equal sample sizes in all treatment by dose combinations, but the experimenter may conduct a study with a sample of $n_i$ for the $i^{th}$ treatment, as long as within each treatment the same number of observations are collected on each dose group or *vise versa*.

**Remark 2.5.** As can be seen from the following counter example, it is encouraging to note that the sufficient condition stated in the above theorem is not necessary. Let $\theta$ be a $2 \times 2$ matrix with $\theta_{i,1} \leq \theta_{i,2}$, $i = 1, 2$, and $\theta_{1,j} \leq \theta_{2,j}$, $j = 1, 2$. Suppose that $\hat{\theta}_{1,1} \leq \hat{\theta}_{1,2}$, $\hat{\theta}_{1,1} \leq \hat{\theta}_{2,1} \leq \hat{\theta}_{2,2}$, and $\hat{\theta}_{1,2} \geq \hat{\theta}_{2,2}$. In this case, for any non-singular weight matrices $W_R$ and $W_C$, the iterative algorithm converges in one cycle, i.e. one column operation followed by one row operation.

**Remark 2.6.** In some situations, such as in the context of ordinal data, the rows (or columns) of the unrestricted estimator $\hat{\theta}$ naturally satisfy the order restriction. In such cases, under the conditions of Theorem 2.1, we need to apply $\mathfrak{C}^{\mathcal{C}}_{W_C}$ on $\hat{\theta}$ (or $\mathfrak{R}^{\mathcal{R}}_{W_R}$) only once.

## 2.3. A simulation study

For an $I \times J$ matrix $\hat{\theta}$ whose components are independently normally distributed, we compared the performance of the proposed order-restricted matrix estimator $\hat{\hat{\theta}}$ with the UMLE $\hat{\theta}$ using a simulation study. Although we considered a variety of order restrictions, patterns of means, sample sizes and weight matrices, in this article we provide the results of a small subset since very similar results were obtained across all patterns.

In the simulation study reported here, $E(\hat{\theta}) = \theta$ with $\theta_{i,j} = i + j$, $1 \leq i \leq I$, $1 \leq j < J - 2$, and $\theta_{i,j} = i - j + J + 1$, $1 \leq i \leq I$, $j \geq J - 2$. Thus we have a simple order along the columns of $\theta$, with $\theta_{i_1,j} \leq \theta_{i_2,j}$, for all $1 \leq i_1 \leq i_2 \leq I$, $j = 1, 2, \ldots, J$, and a tree order along the rows of $\theta$, with $\theta_{i,1} \leq \theta_{i,j}$, for all $1 < j \leq J$ and $1 \leq i \leq I$. Each of the normal variables was generated with a standard deviation of 1 and we chose sample sizes $n_{i,j} = i$, $i = 1, 2, \ldots, I$, and $j = 1, 2, \ldots, J$ so that Variance $(\hat{\theta}_{i,j}) = 1/i$. The $(i, j)^{th}$ element of the weight matrix $W_R$ is given by $\sqrt{i}$ and the weight matrix $W_C = W_R'$.

Our simulation study is based on 10,000 simulation runs and the results are summarized in Table 2. In addition to comparing the average bias of $\hat{\hat{\theta}}$ with that of $\hat{\theta}$, we also computed the percentage reduction in quadratic and quartic loss due to $\hat{\hat{\theta}}$. The reduction in loss relative to $\hat{\theta}$ is defined as $100 \times (1 - \frac{\sum_{i=1}^{I} \sum_{j=1}^{J} E(\hat{\hat{\theta}}_{i,j} - \theta_{i,j})^\delta}{\sum_{i=1}^{I} \sum_{j=1}^{J} E(\hat{\theta}_{i,j} - \theta_{i,j})^\delta})$, where $\delta = 2$ corresponds to quadratic loss and $\delta = 4$ corresponds to quartic loss.

Observe that the proposed procedure reduces the average quadratic loss and quartic loss substantially, without costing much in terms of bias.

## 3. Testing hypotheses under order restrictions

In some applications, such as in Wormser et al. [26], researchers are interested in performing tests of hypotheses regarding the elements of each column of $\theta$ when the rows are subject to order restrictions. The order restrictions on the rows are not part of the hypothesis, but they are present due to the underlying probability model or for other reasons. Thus, the hypotheses of interest are

$$H_0 : \theta_{1,j} = \theta_{2,j} = \cdots = \theta_{I,j}, \ 1 \leq j \leq J,$$

TABLE 2
*Bias and reduction in loss due to the proposed estimator $\hat{\hat{\theta}}$ relative to UMLE $\hat{\theta}$*

| $I$ | $J$ | Average bias | | Percentage reduction relative to $\hat{\theta}$ | |
| --- | --- | --- | --- | --- | --- |
| | | $\hat{\hat{\theta}}$ | $\hat{\theta}$ | Quadratic loss | Quartic loss |
| 2 | 5 | −0.0774 | 0.0014 | 29.80 | 53.54 |
| 2 | 10 | −0.0608 | −0.0009 | 22.38 | 41.41 |
| 5 | 5 | −0.0612 | 0.0018 | 28.12 | 55.06 |
| 5 | 10 | −0.0324 | −0.0011 | 22.75 | 45.02 |

$$(3.1) \quad H_a : (\theta_{1,j}, \theta_{2,j}, \ldots, \theta_{I,j})' \in \Lambda_{k=1}^p \mathcal{M}_k, \ 1 \leq j \leq J,$$

where $\mathcal{M}_k$ are maximally linked subgraphs. However, each row of $\theta$ is itself subject to the restriction $(\theta_{i,1}, \theta_{i,2}, \ldots, \theta_{i,J})' \in \Lambda_{k=1}^q \mathcal{N}_k, \ 1 \leq i \leq I$, where $\mathcal{N}_k$ are maximally linked subgraphs.

In other examples researchers may be interested in testing

$$H_0 : \theta_{i,j} = \theta_{i',j'} \ \forall (i,j) \neq (i',j'), 1 \leq i, i' \leq I, 1 \leq j, j' \leq J,$$

$$(3.2) \quad H_a : \theta \in \mathcal{RAC},$$

where $\mathcal{C} = \Lambda_{k=1}^p \mathcal{M}_k$ and $\mathcal{R} = \Lambda_{k=1}^q \mathcal{N}_k$, and each $\mathcal{M}_k$ and $\mathcal{N}_k$ are suitable maximally linked subgraphs.

In both (3.1) and (3.2) the point estimators of $\theta$ are the same, and are obtained under the order restrictions $\mathcal{RAC}$ using the methodology introduced in Section 2, although the test statistics are different.

For a maximally linked subgraph $\mathcal{M}_k$, defined by $\theta_{s,k} \leq \theta_{s+1,k} \cdots \leq \theta_{r,k}$, the two parameters $\theta_{s,k}$ and $\theta_{r,k}$, which are at the ends of the graph, are said to be *farthest linked parameters* of the subgraph. Under this maximally linked subgraph let $\hat{\hat{\theta}}_{s,k}, \hat{\hat{\theta}}_{s+1,k} \ldots, \hat{\hat{\theta}}_{r,k}$ denote the estimated value of the parameter $(\theta_{s,k}, \theta_{s+1,k}, \ldots, \theta_{r,k})'$ using the methodology described in Section 2. To test the hypotheses given in (3.1), for each maximally linked subgraph, compute the estimated difference between the two farthest linked parameters of the subgraph and divide it by the standard error of the difference under the null hypothesis and take the largest over all maximally linked subgraphs. More precisely we propose the following test statistic for testing (3.1):

$$(3.3) \quad T_1 = \max_{\mathcal{M}_k} \frac{\hat{\hat{\theta}}_{r,k} - \hat{\hat{\theta}}_{s,k}}{\hat{se}(\hat{\theta}_{r,k} - \hat{\theta}_{s,k})},$$

where $\theta_{r,k}$ and $\theta_{s,k}$ are the farthest linked parameters in the maximally linked subgraph $\mathcal{M}_k$ and the max is taken over all maximally linked subgraphs $\mathcal{M}_k$.

Similarly, we may test

$$H_0 : \theta_{i,1} = \theta_{i,2} = \cdots = \theta_{i,J}, \ 1 \leq i \leq I,$$

$$(3.4) \quad H_a : (\theta_{i,1}, \theta_{i,2}, \ldots, \theta_{i,J})' \in \Lambda_{k=1}^p \mathcal{N}_k, \ 1 \leq i \leq I,$$

using the test statistic

$$(3.5) \quad T_2 = \max_{\mathcal{N}_k} \frac{\hat{\hat{\theta}}_{k,r} - \hat{\hat{\theta}}_{k,s}}{\hat{se}(\hat{\theta}_{k,r} - \hat{\theta}_{k,s})},$$

where $\theta_{k,r}$ and $\theta_{k,s}$ are the farthest linked parameters in the maximally linked subgraph $\mathcal{N}_k$ and the max is taken over all maximally linked subgraphs $\mathcal{N}_k$. Then the hypothesis on rows and columns (3.2) may be tested using $T = \max\{T_1, T_2\}$. In the above expressions $\hat{se}$ is computed under $H_0$. A suitable such estimate is derived in Section 4 for the case of ordinal data. The critical values and p-values are obtained by the bootstrap methodology as follows. For each bootstrapped dataset, the $i^{th}$ independent group is formed by taking a simple random sample (with replacement), of appropriate size, from the pooled sample of all subjects across all independent groups. In this resampling procedure, we sample the entire record of a given subject. For each bootstrap sample we compute the test statistic $T$, which is denoted by $T^*$. The sampling distribution of $T^*$ is obtained by generating a large number of bootstrap samples, say 10,000. Then the bootstrap p-value is computed as the proportion of times $T^*$ exceeded the observed $T$.

## 4. Simulation study

We conducted a simulation study to investigate the performance of the proposed bootstrap test based on $T_1$ for comparing $I$ treatment groups when the responses are measured on $J$ ordered categories. In a random sample of $n$ observations on the $i^{th}$ treatment, $i = 1, 2, \ldots, I$, let $X_{i,j}$ denote the frequency of responses in the $j^{th}$ ordered category, $j = 1, 2, \ldots, J$. Further, let $E(X_{i,j}) = n\pi_{i,j}$ and let $\theta_{i,j} = \sum_{k=1}^{j} \pi_{i,k}$ denote the cumulative probabilities with $\theta_{i,J} = 1$. Let $\hat{\pi}_{i,j}$ denote the sample proportion $X_{i,j}/n$ and let $\hat{\theta}_{i,j}$ denote the corresponding cumulative sum of sample proportions. Thus under the multinomial model $\hat{\theta}$ is the UMLE of $\theta$. Note that the rows of $\hat{\pi}$, and hence the rows of $\hat{\theta}$, are independently distributed.

We performed simulations to study the size and power of $T_1$ when testing $H_0$ against $H_a - H_0$, where $H_0 : \theta \in \Theta_0 = \{\theta \mid \theta_{r,j} = \theta_{s,j}, 1 \leq r, s \leq I, j = 1, 2, \ldots, J\}$ and $H_a : \theta \in \Theta = \{\theta \mid r \leq s \Rightarrow \theta_{r,j} \leq \theta_{s,j}\}$. Ordinal data was generated for a variety of parameter configurations (Table 3) with sample sizes of 10, 20, and 50 subjects for each independent group.

Under the null hypothesis, we estimate the standard error of $(\hat{\theta}_{r,j} - \hat{\theta}_{s,j})$ required in $T_1$ by $\hat{se}(\hat{\theta}_{r,j} - \hat{\theta}_{s,j}) = \sqrt{2\tilde{V}_j}$, where $\tilde{V}_j = \frac{1}{n}\sum_{r=1}^{j}\sum_{s=1}^{j}[\tilde{\pi}_r(1-\tilde{\pi}_r)I(r=s) - \tilde{\pi}_r\tilde{\pi}_s I(r \neq s)]$, and $\tilde{\pi}_r = \frac{\sum_{i=1}^{I} X_{i,r} + \frac{I\sqrt{n}}{J}}{nI + I\sqrt{n}}$, $r = 1, 2, \ldots, J$, a pooled Bayes estimator under quadratic loss and suitable Dirichlet distribution prior (Lehmann [15], page 293). Since the usual pooled MLE can take a value of 0 or 1 with a positive probability, we prefer the above estimator over the MLE for calculating $\tilde{V}_j$.

For each configuration listed in Table 3, we performed 10,000 simulations to evaluate the size and power of $T_1$. Critical values of $T_1$ were determined using 10,000 bootstrap samples for each simulation run. The nominal size was set at $\alpha = 0.05$. We compared the proposed test with the order-restricted methods of Grove [9] and Nair [16] as well as with the one-sided Kolmogorov-Smirnov test. For simplicity, we bootstrapped the critical values of Kolmogorov-Smirnov test. Since the procedure of Grove [9] is designed for comparing only two groups, we compared the performance our procedure with Grove [9] for $I = 2$ only. Similarly, comparisons with the one-sided Kolmogorov-Smirnov test was limited to the case $I = 2$ only.

Results of the simulation study are summarized graphically in Figure 1 using scatter plots. The top three panels correspond to the type 1 errors for the three different sample sizes per test group, i.e., $n = 10, 20, 50$. The bottom three panels correspond to the power for the corresponding sample sizes. In each panel the vertical axis represents the estimated probability of rejection of null hypothesis by the proposed critical region, while the horizontal axis represents the estimated probability of rejection of null hypothesis by the three alternative procedures. Thus the six scatter plots represent the comparison between the proposed and each of the three competing test procedures. The '+' symbol corresponds to the comparison between the proposed test and Grove's test (Grove [9]), '#' corresponds to comparison between the proposed test and Nair's test (Nair [16]) and '*' corresponds to the comparison between the proposed test and the one-sided Kolmogorov-Smirnov test. In each scatter plot the diagonal line corresponds to the line of equality. Additionally, in the top three panels a horizontal and a vertical line is provided at $0.05 + \sqrt{(.05 \times .95)/10000}$ to indicate points that exceed the nominal value of 0.05 by at least one standard error.

TABLE 3
*Multinomial parameter configurations used in the simulation study*

| | | Group 1 | | | | | Group 2 | | | | | Group 3 | | | |
|---|---|---|---|---|---|---|---|---|---|---|---|---|---|---|---|
| | | $\pi_{11}$ | $\pi_{12}$ | $\pi_{13}$ | $\pi_{14}$ | $\pi_{15}$ | $\pi_{21}$ | $\pi_{22}$ | $\pi_{23}$ | $\pi_{24}$ | $\pi_{25}$ | $\pi_{31}$ | $\pi_{32}$ | $\pi_{33}$ | $\pi_{34}$ | $\pi_{35}$ |
| $H_0$ | $I=2, J=3$ | 0.33 | 0.33 | 0.33 | | | 0.33 | 0.33 | 0.33 | | | | | | | |
| | | 0.10 | 0.40 | 0.50 | | | 0.10 | 0.40 | 0.50 | | | | | | | |
| | | 0.40 | 0.40 | 0.20 | | | 0.40 | 0.40 | 0.20 | | | | | | | |
| | | 0.01 | 0.49 | 0.50 | | | 0.01 | 0.49 | 0.50 | | | | | | | |
| | | 0.49 | 0.01 | 0.50 | | | 0.49 | 0.01 | 0.50 | | | | | | | |
| | | 0.98 | 0.01 | 0.01 | | | 0.98 | 0.01 | 0.01 | | | | | | | |
| | | 0.01 | 0.98 | 0.01 | | | 0.01 | 0.98 | 0.01 | | | | | | | |
| | | 0.01 | 0.01 | 0.98 | | | 0.01 | 0.01 | 0.98 | | | | | | | |
| | $I=2, J=4$ | 0.25 | 0.25 | 0.25 | 0.25 | | 0.25 | 0.25 | 0.25 | 0.25 | | | | | | |
| | | 0.20 | 0.10 | 0.30 | 0.40 | | 0.20 | 0.10 | 0.30 | 0.40 | | | | | | |
| | | 0.97 | 0.01 | 0.01 | 0.01 | | 0.97 | 0.01 | 0.01 | 0.01 | | | | | | |
| | | 0.01 | 0.97 | 0.01 | 0.01 | | 0.01 | 0.97 | 0.01 | 0.01 | | | | | | |
| | | 0.01 | 0.01 | 0.97 | 0.01 | | 0.01 | 0.01 | 0.97 | 0.01 | | | | | | |
| | | 0.01 | 0.01 | 0.01 | 0.97 | | 0.01 | 0.01 | 0.01 | 0.97 | | | | | | |
| | $I=3, J=4$ | 0.25 | 0.25 | 0.25 | 0.25 | | 0.25 | 0.25 | 0.25 | 0.25 | | 0.25 | 0.25 | 0.25 | 0.25 | |
| | | 0.10 | 0.20 | 0.30 | 0.40 | | 0.10 | 0.20 | 0.30 | 0.40 | | 0.10 | 0.20 | 0.30 | 0.40 | |
| | | 0.40 | 0.30 | 0.20 | 0.10 | | 0.40 | 0.30 | 0.20 | 0.10 | | 0.40 | 0.30 | 0.20 | 0.10 | |
| | | 0.01 | 0.49 | 0.49 | 0.01 | | 0.01 | 0.49 | 0.49 | 0.01 | | 0.01 | 0.49 | 0.49 | 0.01 | |
| | | 0.49 | 0.01 | 0.49 | 0.01 | | 0.49 | 0.01 | 0.49 | 0.01 | | 0.49 | 0.01 | 0.49 | 0.01 | |
| | | 0.97 | 0.01 | 0.01 | 0.01 | | 0.97 | 0.01 | 0.01 | 0.01 | | 0.97 | 0.01 | 0.01 | 0.01 | |
| | | 0.01 | 0.97 | 0.01 | 0.01 | | 0.01 | 0.97 | 0.01 | 0.01 | | 0.01 | 0.97 | 0.01 | 0.01 | |
| | | 0.01 | 0.01 | 0.97 | 0.01 | | 0.01 | 0.01 | 0.97 | 0.01 | | 0.01 | 0.01 | 0.97 | 0.01 | |
| | | 0.01 | 0.01 | 0.01 | 097 | | 0.01 | 0.01 | 0.01 | 0.97 | | 0.01 | 0.01 | 0.01 | 0.97 | |
| | $I=3, J=5$ | 0.20 | 0.20 | 0.20 | 0.20 | 0.20 | 0.20 | 0.20 | 0.20 | 0.20 | 0.20 | 0.20 | 0.20 | 0.20 | 0.20 | 0.20 |
| | | 0.01 | 0.48 | 0.01 | 0.49 | 0.01 | 0.01 | 0.48 | 0.01 | 0.49 | 0.01 | 0.01 | 0.48 | 0.01 | 0.49 | 0.01 |
| $H_1 - H_0$ | $I=2, J=3$ | 0.10 | 0.30 | 0.60 | | | 0.15 | 0.35 | 0.50 | | | | | | | |
| | | 0.10 | 0.30 | 0.60 | | | 0.60 | 0.30 | 0.10 | | | | | | | |
| | | 0.40 | 0.40 | 0.20 | | | 0.80 | 0.10 | 0.10 | | | | | | | |
| | | 0.01 | 0.01 | 0.98 | | | 0.20 | 0.20 | 0.60 | | | | | | | |
| | | 0.01 | 0.01 | 0.98 | | | 0.98 | 0.01 | 0.01 | | | | | | | |
| | $I=2, J=4$ | 0.10 | 0.30 | 0.20 | 0.40 | | 0.30 | 0.10 | 0.40 | 0.20 | | | | | | |
| | | 0.10 | 0.10 | 0.10 | 0.70 | | 0.40 | 0.05 | 0.05 | 0.50 | | | | | | |
| | | 0.10 | 0.20 | 0.30 | 0.40 | | 0.15 | 0.25 | 0.30 | 0.30 | | | | | | |
| | | 0.10 | 0.20 | 0.30 | 0.40 | | 0.40 | 0.30 | 0.20 | 0.10 | | | | | | |
| | | 0.25 | 0.25 | 0.25 | 0.25 | | 0.35 | 0.30 | 0.30 | 0.05 | | | | | | |
| | | 0.01 | 0.01 | 0.01 | 0.97 | | 0.97 | 0.01 | 0.01 | 0.01 | | | | | | |
| | $I=3, J=4$ | 0.10 | 0.30 | 0.20 | 0.40 | | 0.20 | 0.20 | 0.30 | 0.30 | | 0.30 | 0.10 | 0.40 | 0.20 | |
| | | 0.10 | 0.10 | 0.10 | 0.70 | | 0.40 | 0.05 | 0.05 | 0.50 | | 0.60 | 0.10 | 0.10 | 0.20 | |
| | | 0.10 | 0.20 | 0.30 | 0.40 | | 0.10 | 0.20 | 0.30 | 0.40 | | 0.40 | 0.30 | 0.20 | 0.10 | |
| | | 0.10 | 0.20 | 0.30 | 0.40 | | 0.20 | 0.20 | 0.30 | 0.30 | | 0.40 | 0.30 | 0.20 | 0.10 | |
| | | 0.10 | 0.20 | 0.30 | 0.40 | | 0.40 | 0.30 | 0.20 | 0.10 | | 0.40 | 0.30 | 0.20 | 0.10 | |
| | I=3, J=5 | 0.10 | 0.20 | 0.30 | 0.30 | 0.10 | 0.10 | 0.30 | 0.30 | 0.20 | 0.10 | 0.20 | 0.30 | 0.30 | 0.10 | 0.10 |
| | | 0.10 | 0.10 | 0.10 | 0.10 | 0.60 | 0.10 | 0.10 | 0.20 | 0.20 | 0.40 | 0.20 | 0.20 | 0.20 | 0.20 | 0.20 |
| | | 0.10 | 0.10 | 0.10 | 0.10 | 0.60 | 0.50 | 0.10 | 0.10 | 0.10 | 0.20 | 0.90 | 0.01 | 0.01 | 0.01 | 0.07 |

We notice that in general all three test procedures approximately maintain the nominal size of 0.05. In situations involving rare events (e.g. – probability vector (0.01, 0.01, 0.98) all tests are conservative, but even in that case, relative to others, the proposed test appears to recover the nominal size more quickly as the sample size increases.

As indicated by the points above the diagonal line in the three bottom panels of Figure 1, the proposed test seems to enjoy higher power than its competitors in almost all situations considered in this simulation study. In particular, even for parameter configurations that are very close to the null, the proposed procedure has higher power than its competitors. Additionally, in these cases, the proposed procedure appears to increase in power with sample size at a faster rate than its competitors.

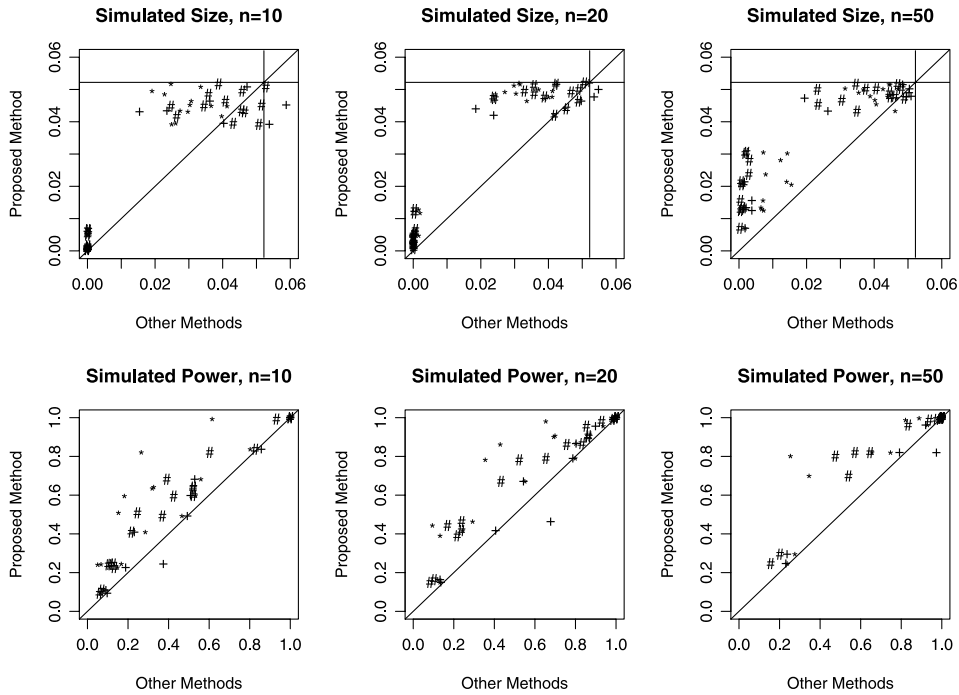

FIG 1. *Power and size comparisons of the proposed procedure with Grove's test (+), Nair's test (#) and Kolmogorov-Smirnov test (*). Nominal size is 0.05. Results for the proposed method are plotted on the vertical axis, and results for other methods are plotted on the horizontal axis.*

In addition to a gain in power, a distinct advantage of the proposed test over likelihood ratio type procedures is the ease of implementation for any arbitrary order restriction on the rows and columns. The procedure described in Grove [9] is limited to two groups. As seen in Wang [25], the likelihood ratio type procedures are very challenging to implement as the number of groups increases.

## 5. Illustration

In the experiment of Wormser et al. [26] mentioned in the introduction of this paper, the effect of mustard gas on the skin of mice was evaluated using 6 ordinal variables, namely, *subepidermal microblister formation, epidermal ulceration, epidermal necrosis, acute inflammation, hemorrhage, and dermal necrosis*. The experiment consisted of exposing 10 mice of each genotype (i.e. COX-2-d, WT and COX-1-d) to 0.317mg of sulfur mustard. Changes in skin condition of each animal (as measured by the above 6 variables) were noted on ordinal scale ranging from "unremarkable", "minimal", "mild", "moderate", to "marked". For each response variable and each genotype, in Table 4 we provide the sample cumulative proportion of animals in each category.

TABLE 4
*Cumulative relative frequencies for each genotype each for response variable*

| | | Level of skin injury | | | | |
|---|---|---|---|---|---|---|
| **Response variable** | **Genotype** | Unremarkable | Minimal | Mild | Moderate | Marked |
| Microblister | COX-1-d | 0 | 0.2 | 1 | 1 | 1 |
| | WT | 0 | 0.5 | 0.9 | 1 | 1 |
| | COX-2-d | 0.3 | 0.8 | 1 | 1 | 1 |
| Ulceration | COX-1-d | 0.4 | 0.5 | 0.8 | 0.9 | 1 |
| | WT | 0.8 | 0.9 | 1 | 1 | 1 |
| | COX-2-d | 1 | 1 | 1 | 1 | 1 |
| Epidermal necrosis | COX-1-d | 0 | 0 | 0 | 0.3 | 1 |
| | WT | 0 | 0 | 0.1 | 0.8 | 1 |
| | COX-2-d | 0.1 | 0.2 | 0.7 | 0.8 | 1 |
| Acute inflammation | COX-1-d | 0 | 0 | 0.5 | 1 | 1 |
| | WT | 0 | 0 | 0.6 | 1 | 1 |
| | COX-2-d | 0.1 | 0.5 | 1 | 1 | 1 |
| Hemorrhage | COX-1-d | 0 | 0.1 | 0.9 | 1 | 1 |
| | WT | 0 | 0 | 1 | 1 | 1 |
| | COX-2-d | 0.1 | 0.4 | 1 | 1 | 1 |
| Dermal necrosis | COX-1-d | 0 | 0.1 | 0.8 | 1 | 1 |
| | WT | 0 | 0.1 | 0.9 | 1 | 1 |
| | COX-2-d | 0.2 | 0.4 | 1 | 1 | 1 |

For each response variable, the statistical hypothesis of interest was motivated by the underlying biology. COX-2 is involved in a variety of inflammatory processes caused by noxious stimuli. For instance, COX-2 is induced within 12 hours and persists up to 3 days after excisional injury in rat skin. COX-2 expression is seen in the basal cell layer, peripheral cells in the outer root sheath of hair follicles, and in fibroblast-like cells and capillaries near the wound edges Lee et al. [13].

Neutrophil COX-2 protein expression after burn-induced injury in mice significantly increased at 4 hours and dramatically decreased 36 hours after injury (He et al. [10]). Due to the central role of COX-2 in the inflammatory processes, it is believed that the effect of mustard gas on COX-2 deficient animals will tend to be less severe than that on the Wild Type animals.

Nevertheless, COX-1 may have a protective effect on the skin against sulfur mustard as shown in other organs like the kidney and brain (Lin et al. [14]) and Vane et al. [24]. COX-1 was suggested to confer protection on the epithelial cells of the crypts of Lieberkühn, through promotion of crypt stem cell survival and proliferation, in the ileum of irradiated mice Cohn et al. [5].

As in our skin model case, Cohn et al. [5] also concluded that prostaglandins produced through the COX-1 pathway, may not be important in unstressed conditions, but still may have a protective role in the response to epithelial injury. Therefore, we expect mice with COX-1 deficiency to have a more severe response, on average, than wild-type animals. Then, the experimental setting is the comparison of three groups with six responses, each of which is measured on an ordinal scale with five levels and subject to a simple order.

We analyzed each response variable separately but adjusted the p-value using Bonferroni correction. The analysis of each response variable was carried out using the methodology developed in this paper. Based on the above discussion, for each of the 6 response variables, we tested the following hypotheses, where $i = 1$ represents COX-1-d, $i = 2$ represents WT, and $i = 3$ represents COX-2-d:

$$H_0 : \theta_{1,j} = \theta_{2,j} = \theta_{3,j}, \ j = 1, \ldots, 4.$$
$$H_a : \theta_{1,j} \leq \theta_{2,j} \leq \theta_{3,j}, \ j = 1, \ldots, 4.$$

We computed p-values by bootstrapping the test statistic with 50000 resamplings, each time sampling the entire record of an animal to preserve the underlying dependence structure between the 6 response variables.

After performing Bonferroni correction for 6 tests, significant genotype trends were found in subepidermal microblister formation ($p = 0.0310$), ulceration ($p = 0.0077$), epidermal necrosis ($p = 0.0034$), and acute inflammation ($p = 0.0247$). We failed to see significant trends in hemorrhage ($p = 0.1181$) and in dermal necrosis ($p = 0.5170$).

## 6. Discussion

In this article we have extended the iterative algorithm of Dykstra and Robertson [7] to arbitrary order restrictions on the rows and columns of a matrix as long as each row is subject to the same order restriction and each column is subject to the same order restriction. However, the order restrictions on the rows need not be same as those on the columns. Within each row/column the new algorithm makes use of the same estimation procedure introduced in Hwang and Peddada [11]. If rows and columns are subject to a simple order then, for independently and normally distributed data, the proposed algorithm is identical to the algorithm of Dykstra and Robertson [7]. We derive a sufficient condition for the algorithm to converge in a single application of Hwang and Peddada [11] method on rows and on columns. As an example, the sufficient condition is satisfied in a balanced design.

Using the point estimators derived in this paper, we introduced a new test statistic which is a Kolmogorov type distance on the graph of order restrictions as in (Peddada et al. [17]). The new procedure is easy to implement and simulations performed in this paper suggest that it has higher power in most cases than some of the existing procedures. A part of the reason for the new procedure to perform better than some of the standard procedures, such as the Kolmogorov-Smirnov test, is because it uses improved order-restricted point estimators. As seen from our simulation studies, the new estimator enjoys substantially smaller quadratic (and quartic) risk than the unrestricted estimator, which is used in the Kolmogorov-Smirnov test.

An added advantage of the proposed methodology is that it is applicable to a very broad collection of matrix order restrictions and is computationally easy to implement.

**Acknowledgments.** The authors wish to thank David Dunson, Grace Kissling, the reviewers and the editor for providing useful comments that greatly improved the presentation of the manuscript.